# Leveraging guidelines for ethical practice of statistics and computing to develop standards for ethical mathematical practice: A White Paper

Catherine Buell[1]*, Victor I. Piercey[2]*, Rochelle E. Tractenberg[3]*

* Authors contributed equally to this manuscript and names appear in alphabetical order.

1 Department of Mathematics, Fitchburg State University, 60 Pearl Street, Fitchburg, MA USA 01420-2697
2 Department of Mathematics, Ferris State University, 1201 S. State Street, Big Rapids, Michigan USA 49307
3 Departments of Neurology; Biostatistics, Bioinformatics & Biomathematics; & Rehabilitation Medicine, Georgetown University, Suite 207 Building D, 4000 Reservoir Road NW, Washington, DC USA 20057.

**Running head:** Ethical Mathematical Practice

This White Paper describes the incubation project funded by National Science Foundation to the co-authors (Proposal Numbers 2024227 -VIP; 2024279-RET; and 2024244-CB)

Cite this White paper: Buell C, Piercey VI, Tractenberg RE. (2022). *Leveraging guidelines for ethical practice of statistics and computing to develop standards for ethical mathematical practice: A White Paper*, ArXiv: 2209.09311 [math.HO] https://doi.org/10.48550/arXiv.2209.09311

Final version (March 2024) of the Proto Ethical Guidelines are available as part of the OSF project: https://osf.io/ahj9q/

**NOTE: This White Paper fully describes the work of the NSF grant. The shorter manuscript is forthcoming in Journal of Science & Engineering Ethics (2024). That manuscript overlaps significantly with this one but features a different introduction and contains the final version of the Ethical Guidelines described in this white paper. That manuscript is available as a preprint:** Tractenberg RE, Piercey VI, Buell C. (2024). Defining “ethical mathematical practice” through engagement with discipline-adjacent practice standards and the mathematical community. SocArXiv: https://doi.org/10.31235/osf.io/6a5f7

**ABSTRACT**

We report the results of our NSF-funded project in which we alpha- and beta-tested a survey comprising all aspects of the ethical practice standards from two disciplines with relevance to mathematics, the American Statistical Association (ASA) and Association of Computing Machinery (ACM). Items were modified so that text such as “A computing professional should…” became, “The ethical mathematics practitioner…”. We also removed elements that were duplicates or were deemed unlikely to be considered relevant to mathematical practice even after modification. Starting with more than 100 items, plus 10 demographic questions, the final survey included 52 items (plus demographics), and 142 individuals responded to our invitations (through listservs and other widespread emails and announcements) to participate in this 30-minute survey. This white paper reports the project methods and findings regarding the community perspective on the 52 items, specifically, which rise to the level of ethical obligations, which do not meet this level, and what is missing from this list of elements of ethical mathematical practice. The results suggest that the community of mathematicians perceives a much wider range of behaviors to be subject to ethical practice standards than is currently represented.

**Introduction**

In discussion around "ethics in mathematics", common talking points include recent scandals and international interdisciplinary conversations around ethical and unethical algorithms (see, e.g., O'Neil 2016). These conversations forced some mathematicians to consider what “ethical practice” looks like in their profession (Karaali 2019). Others, like Müller (2018), point out that

> [f]or many pure mathematicians, their research does not entail a normative component even though it changes every time their research is used in the real world. This suggests an increasing gap between how mathematicians think about their field and how society does, or at least how social scientists do, a gap which is deepening over time due to the communicative power of numbers (pp. 5).

But beyond the mathematicians' identity through their research, there is a growing and concerned interest in the ethical practice of mathematics and its culture. Rittberg et al. (2020) analyzed mathematical practice through virtue and epistemic lenses and provided case studies for ethical practice in academia and education considering access, validity, norms, and values. Continuing in education, Ernest (Ernest 2018) speaks to the ethical obligations of how we teach mathematics. Others (eg. Franklin 2005, Ferrini-Mundy 2008) build teaching ethics and research ethics into graduate-level mathematics training. Research/work, community membership, and education are fundamental pieces of mathematical practice and each have ethical considerations for all practitioners.

We start from an assumption that mathematicians *want to be ethical* even when they don't think ethics apply to them. What is needed is what Chiodo and Bursill-Hall (Chiodo and Bursill-Hall 2018) refer to as an "ethical consciousness" among mathematicians - an instinct to ask oneself the ethical and social implications of one's work regardless of whether the applications are obvious. They speak to many ethical obligations, some of which include considering future applications of mathematics research and ethical practices for research funding. Many examples of ethical problems examined in mathematical sciences are related to computer science, algorithms, and applications. Mathematics is at the heart of each of these areas, and each has been aided by abstract problems solved by mathematicians where applications were originally unlikely to occur or may be unseen. For example, number theory is one of the most venerable and abstract areas of mathematical research, but it is critical (i.e., applied) in cryptography. Bass (2006) articulated that engagement, particularly with research/scholarship in mathematics, comprises both a professional and a disciplinary aspect (p. 103-104).

Similarly, "The entire community of scientists and engineers benefits from diverse, ongoing options to engage in conversations about the ethical dimensions of research and (practice)," (Kalichman 2013: 13). However, we (Buell & Piercey 2019) join the growing voice of mathematicians cited above and argue that one must go beyond the implications or applications, and focus on an "ethical consciousness" in the practice of mathematics and the culture of mathematics. Mathematics as a practice goes beyond solely content, including knowledge-production, research, knowledge-sharing, culture, and teaching, which are also featured by Bass (2006) and Rios et al.

(2019). Each has ethical components, some unique to mathematics, but all mathematics practitioners could potentially place themselves in multiple aspects of the ethical framework of mathematics. In short, ethics is an important conversation and all too often mathematicians have been left out (either by choice or by distinction or training). Educators have been part of the conversation about *ethically teaching* mathematics (see Neyland 2004; Neyland 2008; Sowder 1998; Atweh et al. 2012); however, researchers and mathematics educators have only scratched the surface and have not met on sufficient common ground needed address the issues of ethics in mathematical culture. More tools, knowledge, and community-commitment is needed to facilitate the integration of necessary ethical discussions into our classrooms, societies, research, and profession.

**Background**

Mathematicians can either start from scratch to derive a new set of ethical practice standards, or, they can join the discussions *in progress* within statistics and computing, two fields intimately – and already - involved in the ethical use of quantitation and data (leveraging mathematics). Specifically, mathematicians might leverage guidelines for ethical practice of these two related professions. Although statistics and computing rely on foundational mathematics, practitioners in, and users of, each discipline have specific aspects that are unique (see, e.g., Tractenberg 2019). Thus, rather than concluding that mathematics has no content that could be subject to ethical practice standards beyond ethical scholarship and disciplinary preparation (AMS 2019), or beyond avoiding/managing conflicts of interest (MAA 2017b), this study sought to explore the perceptions by the mathematical community of the ethical practice standards long maintained in computing by the Association of Computing Machinery (ACM since 1992) and in statistics by the American Statistical Association (ASA since 1995). We also sought community input on additional ethical considerations apart from what computing and statistics practitioners and users have articulated, to ensure a description of “ethical mathematical practice” that is authentic as well as comprehensive.

Scholars have debated the efficacy of ethics codes (see, e.g., Beauchamp and Bowie 1979, Hoffman et. al. 1984, Weller 1988; see also McNamara et al. 2018), but before determining that codes do not work to promote ethical practice (e.g, McNamara et al. 2018; see also May & Luth 2013; Antes et al. 2010), more and focused efforts are needed to teach and give practice with the use and utility of those codes (Tractenberg, et al. 2015). A code – or set of ethical practice standards - requires a profession to articulate duties and responsibilities of a member of the profession that go beyond compliance with the law (Weller 1988). The presence of a code is evidence of a profession’s attention to their own professional responsibilities. Such a code for mathematics practice could be used to initiate and support the development of “ethical consciousness” among established practitioners as well as those in training.

There are a wide range of member-societies with the goal of advancing mathematical sciences (see the Conference Board of the Mathematical Sciences member organizations). We focused on the American Mathematical Society (the AMS) and the Mathematical Association of America (the MAA) as they are the organizations with the broadest reach in the United States and they have each adopted some form of ethical principles. The AMS serves primarily mathematicians engaged in research as well as business and industry, while the MAA serves mathematicians who engage in both teaching roles and the scholarship of teaching and learning in higher education. Both societies have statements of ethics largely surrounding issues of plagiarism and publication (AMS) and conduct at meetings (MAA). While other issues are addressed, neither the AMS nor the MAA describe ethical practice standards for mathematics generally or comprehensively. This is perhaps not surprising since the significance of ethics in mathematics is new to the community while these organizations are longstanding.

The MAA also has a detailed whistleblower policy (MAA 2017a) as well as a conflict of interest policy for those in official positions (MAA 2017b), but in 2017, the MAA “Code of Ethics” consisted of a single paragraph. This paragraph refers to “ethical business and professional practices” but these practices are not specified, nor is a definition of what constitutes “ethical practices” offered. The MAA text also referred to legal compliance, but ethics and ethical codes are designed to fill in gaps between responsible conduct and the law and are often adopted to allow for a profession to regulate themselves (see, e.g., Weller 1988). In 2021, the MAA revised its "Code of Conduct" for meetings and interactions (MAA 2021). The MAA’s Committee on Undergraduate Programs in Mathematics (CUPM) Guide to Majors in the Mathematical Sciences does not include any references to ethics or responsible conduct (CUPM, 2015). Moreover, neither the whistleblower nor conflicts of interest policies is part of the broader ethics statement, and thus may easily be perceived as *not* being part of "professional practices" as articulated *within* the "ethics statement" of the MAA.

The AMS has a policy related to the treatment of students and early-career employees (AMS 2007) as well as an ethics policy (AMS 2019). Most of the ethics policy is related to plagiarism and publication processes as well as institutional responsibilities in conferring degrees. Section II addresses social responsibility, which include statements against discrimination, addresses conflicts of interest, includes a statement expressing a preference for open mathematics over secrecy with the exception of security concerns, and a statement related to temporary employees.

Given the variability across mathematical societies in their statements of ethical practice, and their specificity to research activities in spite of the very broad applicability of mathematics practice outside of research and academia, we explored guidance from allied disciplines, statistics in the American Statistical Association (ASA) Ethical Guidelines for Statistical Practice (ASA 2018; updated in 2022) and computing through the Association for Computing Machinery (ACM) Code of Ethics (ACM 2018). The ethical practice standards from these two organizations are intended to offer guidance for all practitioners, not only those in academia or those doing research.

This study was initiated to gather input from mathematics practitioners across career stages (early/middle/late), practice setting (academia/government/industry) and role (instructor, scientist, or both) in the mathematical community through focused surveys. Community input was sought to address the following questions about "ethical mathematical practice":

1. Which elements of the existing ethical codes of the AMS and MAA are perceived to be relevant to ethical practice by the mathematics community?

2. What ethical guidelines from math-adjacent professional societies (ACM, ASA) does the mathematics community believe are relevant?

3. What other guidelines are necessary that are unique to mathematics? What ethical guidance is lacking from the AMS, MAA, ACM, and ASA guidelines?

## Methods

This project was granted IRB exemptions from all three participating institutions (Ferris #IRB-FY19-20-205, Fitchburg IRB #202021-14, and Georgetown IRB ID #00002454).

### Item Selection

The project's engagement with the mathematics community began with the creation of a preliminary set of items to be considered as part of a "Proto-Ethical Mathematics Practice Guideline" document. We began with all 52 items from the 2018 ASA Guidelines and the 24 ACM (2018) items. A thematic analysis of the 29 AMS (2019) and seven MAA (2017) codes was carried out to derive the most concise list of items that reflected the most general aspects (i.e., not just research) of ethical mathematical practice. Both the AMS and MAA codes are narrative; so when a theme from either was identified as an actual numbered item in either the ASA or ACM codes, we utilized those numbered items.

The *stems* of items in ACM and ASA practice standards differ: ACM Code of Ethics (2018) items have the stem, “A computing professional should…” while the ASA Ethical Guidelines (2018) items have the stem, “The ethical statistician…”. AMS and MAA content do not include stems.

The thematic analyses of the 110 items across the four ethical code documents yielded a preliminary sample of 86 items to be alpha and beta tested. *Alpha testing* occurred at the virtual Joint Mathematics Meeting where attendees spent 1.5 hours considering subsets of the 86 elements. Fifty people joined the virtual Town Hall meeting where we separated them into six groups. Each group was assigned between 13-16 items from each of these source Guideline documents. Groups went through their lists and indicated whether (yes/no) that item would be considered relevant for “ethical mathematical practice”.

Table 1 presents an elemental list derived from the (narrative) AMS code of ethics (2019).

Table 1. AMS Code of Ethics themes (breakdown of Code into items)

| **I. Mathematical research and its presentation** |
|---|
| ● Do not plagiarize, correct attribution when appropriate is essential. |
| ● Be knowledgeable in your field |
| ● Give appropriate credit |
| ● Do not claim a result in advance of its having been achieved; publish full details of results without unreasonable delay after announcing results. |
| ● Use no language that suppresses or improperly detracts from the work of others |
| ● Correct or withdraw work that is erroneous |
| ● A claim of independence may not be based on ignorance of widely disseminated results |
| ● Ensure appropriate authorship |
| **II. Social responsibility of Mathematicians** |
| ● Encourage and promote mathematical ability without bias and review programs to ensure consideration of a full range of students. |
| ● Avoid conflicts of interest and bias in reviewing, refereeing, or funding decisions. |
| ● Respect referee anonymity |
| ● Resist excessive secrecy, promote dissemination/publication |
| ● Disclose implications of work to employers and the public when work may affect public health, safety, or general welfare. |
| ● Do not exploit workers with temporary employment at low pay/excessive work) |
| **III. Education and Granting of Degrees** |
| ● Granting a degree means certifying competence for work. |
| ● PhD level work is ensured by the degree grantors to be high level and original |
| ● PhD is only awarded to those with sufficient knowledge outside the thesis area. |
| ● Degree grantors must honestly inform degree earners about job market/ employment prospects. |
| **IV. Publications** |
| ● Editors should be reasonably sure of the correctness of articles they accept. |
| ● Editors should ensure timely and current reviews. |
| ● Submissions for review are treated as privileged information. |
| ● Editors must prioritize the first submitted version of a paper. |
| ● Editors must inform authors if there is a delay in potential publication |
| ● Publication cannot be delayed for any reason except the authors' interest/actions. |
| ● Date of submission and revisions must be published with any article. |
| ● Editors must be given/accept full responsibility for their journals, resist outside agency pressures and notify the public of such pressure. |
| ● Editors and referees must respect the confidentiality of all submitted materials as appropriate. |
| ● Mathematical publishers must respect the mathematical community and disseminate work accordingly. |
| ● The American Mathematical Society will not play a role/endorse any research journal where any acceptance criterion conflicts with the principles of the AMS guidelines. |
| Note: Grayed out items were excluded from the survey. Other items were excluded if there were more concrete versions on other lists. |

The Ethical Guidelines of the AMS (2019) were reviewed, and 13 elements (of 29 items abstracted from the Guidelines document) were retained. The 16 items that are grayed out in Table 1 were omitted from the survey by the co-authors for one of two reasons. Most typically they were highly limited to very few mathematics practitioners (e.g., “Date of submission and revisions must be published with any article” or “Editors should ensure timely and current reviews”), so that they would be unlikely to be found relevant to ethical mathematical practice by respondents. An alignment table determined whether the remaining 13 items were also reflected in the ASA Ethical Guidelines (2018, comprising 52 items, plus a preamble and the general statements of the eight Principle domains). Ten of the 13 themes from the AMS were also reflected in unique items on the ASA Ethical Guidelines for Statistical Practice (2018). Any ASA or ACM elements that are endorsed in the survey reflecting the AMS items would be interpreted as endorsement of those AMS items. Three unique AMS items (non-exploitation of workers; honest information about job prospects; certification of quality of Ph.D) were retained, and the definition of plagiarism from the AMS was also retained for the survey (instead of the more general descriptions of the ASA and ACM).

Table 2. MAA Code item-level breakdown

| **I. MAA Code of Conduct** |
| --- |
| ● The MAA is committed to adhering to ethical business and professional practices, and to following a policy of honesty and integrity, in the full range of MAA activities. |
| All employees of the MAA and all members engaging in the business, operations, and activities of the MAA shall adhere to all federal, state, and local laws and regulations and conduct themselves in a proper ethical manner. |
| ● The MAA requires Directors, Officers, Members ,those compensated by the MAA and those donating their time, and all employees to observe high standards of business and personal ethics in the conduct of their duties and responsibilities. |
| ● All employees and representatives of the MAA must practice honesty and integrity in fulfilling their responsibilities and comply with all applicable laws and regulations. |
| **II. MAA Whistleblower Policy** |
| ● The MAA will not tolerate intimidation, coercion, or discrimination of any kind against employees or other individuals who file complaints or who |

| | |
|---|---|
| | testify, assist, or participate in any manner in an investigation or hearing. |
| ● | It is the responsibility of all Directors, Officers, members and employees to comply with the Code of Ethics and to report violations or suspected violations in accordance with this Whistleblower Policy. |
| ● | No Director, Officer, member, or employee who in good faith reports a violation of the Code of Ethics shall suffer harassment, retaliation or adverse employment consequences. |
| **III. MAA Welcoming Environment** | |
| ● | The MAA encourages the free expression and exchange of ideas in an atmosphere of mutual respect and collegiality. No discriminatory, harassing, or threatening by any staff member to any other person engaged in MAA operations or activities. |

Note: Grayed out items were excluded from the survey if they applied to only a small subset of mathematics practitioners, or, if there were more concrete versions on other lists.

The MAA Code of Ethics (2017) was reviewed by the co-authors and deemed not sufficiently specific to ethical practice of mathematics to include in the survey. The ASA and ACM include more explicit items discouraging discrimination, harassment, and other illegal acts – against any person, not only against other members, so those items were included in the survey instead of the MAA items so that their themes were represented in the beta version.

*Beta testing* was accomplished when the authors reduced the starting number of survey items from 86 and any items that Joint Mathematics Meetings (JMM) 2021 Town Hall meeting attendees identified as important, but not already in the list of 86, down to a set of all those items that could be framed with respect to mathematics practice (i.e., by changing “computing professional” to “mathematics practitioner”, or by changing specifically statistics or computing terminology to be more consistent with math instead). A small subset of the JMM 2021 Town Hall attendees agreed to be contacted for input on this version of the survey, to offer input (on the content). These beta testers did not respond to survey items, only reviewed them. As an example of modification, ACM Principle 1.1, “A computing professional should… Contribute to society and to human well-being, acknowledging that all people are stakeholders in computing”, became, “The ethical mathematics practitioner… Recognizes that if they engage in mathematics practice, they do so in a social and cultural context. Practitioners should contribute to society and to human well-being, acknowledging that all people are stakeholders in mathematics.” Beta testers ensured the link to the survey worked, commented on clarifications to the instructions, and also identified typos and other irregularities the investigators missed.

*Final Survey*: At the time of the survey (2021), the ASA was revising its ethical guidelines (Tractenberg et al. 2021). A town hall session (1.5 hours) was approved by the Joint Mathematics Meeting organizers for the virtual 2021 meeting (January 2021). Given input from the JMM 2021 town hall meeting, if items were identified there that had been formulated for inclusion in the ASA revisions, we utilized those new items in the Beta test.

We winnowed the beta list down to 52 total items (plus demographics) after the removal of duplicates, and the determination of which of the beta version items were **unlikely** (in the authors’ or beta-testers’ opinions) to be viewed as relevant to the ethical practice of mathematics. Examples of items omitted include, “(the ethical statistician) Employs selection or sampling methods and analytic approaches appropriate and valid for the specific question to be addressed, so that results extend beyond the sample to a population relevant to the objectives with minimal error under reasonable assumptions.”

(ASA Ethical Guideline (2018) Principle A2); and "(A computing professional should) Access computing and communication resources only when authorized or when compelled by the public good." (ACM Code of Ethics 2.8). Some themes represented across several codes were combined and presented as an item from a single code. Table A in the Appendix denotes the origins of common themes.

The 52 item survey (28 ASA, 14 ACM, 4 Town Hall, 4 AMS, and 2 hybrid ASA/ACM; ASA/Town Hall) was deployed, with permissions, by sharing the survey invitation and URL link to it on SurveyMonkey through online messaging boards for members of several professional organizations including the AMS, MAA (including the Business, Industry, and Government Special Interest Group of the MAA), the Society for Industrial and Applied Mathematics (SIAM), and the American Mathematical Association of Two Year Colleges (AMATYC).

Prefacing the survey was the following statement:

> *The items are derived from several sources, so there is a bit of redundancy, but generally speaking, the items can be grouped as reflecting diverse elements of **mathematical practice**. We define the practice of mathematics to include mathematical work; the context in which or for which the work is done; the role of the practitioner; and the matter to which the mathematical work is directed or applied. The survey asks you to consider whether each of the following items is relevant to the practice of mathematics.*
>
> *Answer YES if you feel the item is an ethical obligation for the ethical mathematics practitioner. Answer NO if you feel the item is relevant, but not an ethical obligation; OR, if you feel the item is not relevant to ethical mathematical* *practice. We have included an option for you to comment on your answer. Be sure to consider yourself as a mathematics practitioner, but also other practitioners in the mathematical community who may have different roles than you.*

Also, all items in the survey had the same stem, "The ethical mathematics practitioner:", for example, the first item would be read as:

**"The ethical mathematics practitioner:**

***1. Works in a manner intended to produce valid, interpretable, and reproducible results."***

To increase interpretability of the survey results, we formulated the survey questions to include just "yes" and "no" answers (rather than a Likert scale of respondent-perceived relevance for each item), asking individuals to simply state whether or not they believe each item (given the stem, "the ethical mathematics practitioner") was or was not "relevant to the practice of mathematics." Each item also included the opportunity to comment on either the item or the participant's response. One final item, "is there anything missing?" was also included.

With the exception of demographics, responses of "Yes" represent agreement/endorsement of an item: i.e., "YES if you feel the item is an ethical obligation for the ethical mathematics practitioner." We contemplated how best to present the endorsement rates by item, including a simple tabulation (see Table A, Appendix) and grouping items by endorsement level in order to better understand community thoughts about the relevance of each item for a new set of guidelines for "ethical mathematical practice". Ultimately we utilized simple cluster analysis (in R) based on the endorsement values, and also studied the comments and missingness patterns in the data, but determined the final grouping based on our shared understanding of "general consensus" to be represented by 85%-100% endorsement. The lowest levels of agreement would comprise 0-69.9% endorsement, and thus our middle range of endorsement ended up being 70.0%-84.9%.

**Item-level analysis**

Summary of survey items was based on the numbers of respondents who endorsed each item, but we also conducted an informal thematic analysis of any comments that were obtained on each item, as well as on the "is anything missing?" responses.

**Results**

**General and demographic results**

The survey was open for responses for three months. A total of 142 individuals completed the survey. Their demographics are presented in Figure 1 and Table 3 (subtables A-E) below.

Figure 1. Sample demographics

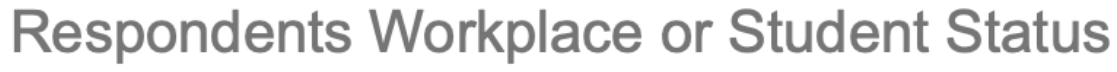

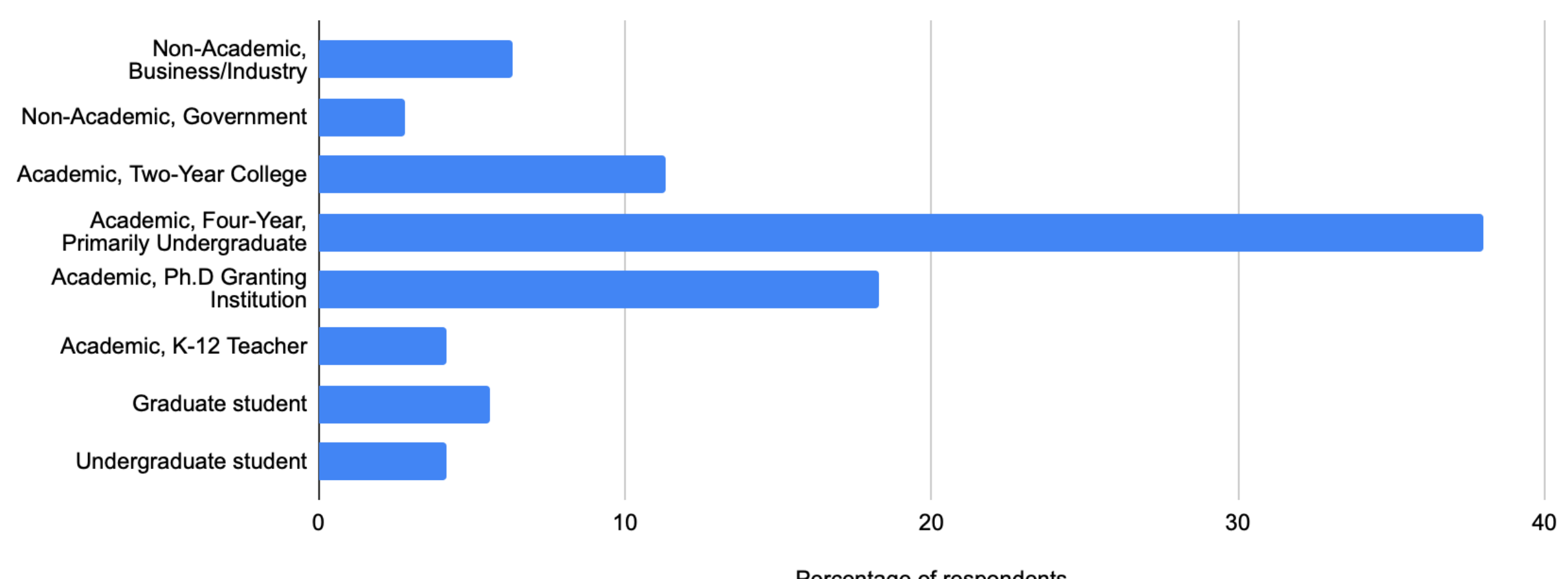

Table 3. **Sub-tables A-E describing respondent demographics**

| **A. Workplace or student status** | Frequency | % |
|---|---|---|
| Non-Academic, Business/Industry | 9 | 6.3% |
| Non-Academic, Government | 4 | 2.8 |
| Academic, Two-Year College Post-Doc, Faculty, or Administration | 16 | 11.3 |
| Academic, Four-Year, Primarily Undergraduate College/University Post-Doc, Faculty, or Administration | 54 | 38.0 |
| Academic, Ph.D Granting Institution Post-Doc, Faculty, or Administration | 26 | 18.3 |
| Academic, K-12 Teacher | 6 | 4.2 |
| Graduate student | 8 | 5.6 |
| Undergraduate student | 6 | 4.2 |

| **B. Highest degree** | Frequency | % |
|---|---|---|
| Bachelors | 4 | 2.8 |
| Ed.D | 1 | 0.7 |
| Masters | 9 | 6.3 |
| Ph.D./Ph.D in progress | 82 | 57.7 |
| Other | 3 | 0.7 |
| Missing | 43 | 30.3 |

| **C. Years of experience** | Frequency | % |
|---|---|---|
| 0 - 5 | 17 | 12.0 |
| 6 - 10 | 10 | 7.0 |
| 11 - 15 | 7 | 4.9 |
| 16 - 20 | 12 | 8.5 |
| 21 - 25 | 13 | 9.2 |
| 26 - 30 | 8 | 5.6 |
| More than 30 | 31 | 21.8 |
| Missing | 44 | 31.0 |

| **D. Gender** | Frequency | % |
|---|---|---|
| Man | 57 | 40.1 |
| Woman | 33 | 23.2 |
| Non-binary | 3 | 2.1 |
| Prefer not to disclose | 5 | 3.5 |
| Missing | 44 | 31.0 |

| **E. Ethnicity** | Frequency | % |
|---|---|---|
| Asian or Asian American | 3 | 2.1 |
| Black or African American | 3 | 2.1 |
| Hispanic or Latino | 4 | 2.8 |
| White or Caucasian | 78 | 54.9 |
| Prefer not to respond | 5 | 3.5 |
| Missing | 49 | 34.5 |

**Endorsement rates**

Of the 52 items - one had 100% endorsement: “Discloses conflicts of interest, financial and otherwise, and manages or resolves them according to established (institutional/regional/local) rules and laws.” This item wording came from the ASA Ethical Guidelines, but is also mentioned on the codes of the AMS, and ACM (albeit less explicitly), and was the only item with unanimous responses. Several others had endorsement at 99%, including “Accepts full responsibility for their own work; does not take credit for the work of others; and gives credit to those who contribute. Respects and acknowledges the intellectual property of others" and "Avoids plagiarism".

Of the 52 items, 50% (26/52) were endorsed by 85-100% of respondents, with 38% (20 items) being endorsed by 90-100% of respondents. A further 17 items (36.7%) were endorsed by 70-84.9% of the sample. Thus, 43/52 or 82.7% of the items on the survey were endorsed by at least 70% of respondents as being "an ethical obligation for the ethical mathematics practitioner". All of the 13 items reflecting AMS content, plus two of four items suggested from the Town Hall, were perceived to be relevant to ethical practice by at least 70% of respondents from the mathematics community. Three items unique to the ethical practice of statistics (ASA), and four unique to ethical computing (ACM), plus two suggested by mathematics community members at the Town Hall, were endorsed by 38.8-69.9% of respondents.

Tables 4, 5, and 6 below represent items endorsed at above 85%, between 70% and 84.9%, and below 69.9% (respectively). This grouping is based on our understanding of "general consensus" as neither an analysis of the item-by-item thematic comments nor the use of a clustering algorithm provided splitting points that used intrinsic qualities of the data. In the tables below, an element is labeled by color to denote its primary origin: purple from the AMS code, peach from the ACM code, blue from the JMM 2021 meeting, and green from the ASA code. Table A in the Appendix provides a more detailed picture of origin. Additionally, as the numbered items increase towards 52, the number of missing responses increases up to 44 (ie, 44 responses were missing at items 51 & 52).  The number of missing responses is presented in parenthesis next to the "% agree" as the number changes.

Table 4 features the endorsement results for the 26 items with the highest levels of agreement (85%-100%).

**Table 4: Survey items with 85-100% agreement**

| | The ethical mathematics practitioner... | % agree (# not responding) |
|---|---|---|
| 1 | Works in a manner intended to produce valid, interpretable, and reproducible results. | 94.3 (1) |
| 2 | Does not knowingly accept work for which they are not sufficiently qualified, is honest about any limitation of expertise, and consults others when necessary or in doubt. It is essential that mathematics practitioners treat others with respect. | 87.9 |
| 3 | Identifies and mitigates any efforts to predetermine or influence the results or outcomes of mathematical practices; resists pressure to solve unethical problems/support predetermined outcomes. | 94.3 |
| 4 | Accepts full responsibility for their own work; does not take credit for the work of others; and gives credit to those who contribute. Respects and acknowledges the intellectual property of others. | 99.3 |
| 6 | Discloses conflicts of interest, financial and otherwise, and manages or resolves them according to established (institutional/regional/local) rules and laws. | 100 |
| 7 | Is candid about any known or suspected limitations, assumptions, or biases when working with data, methods, or models. Objective and valid interpretation of the results requires that the underlying analysis recognizes and acknowledges the degree of reliability and integrity of the data or the model. | 97.9 |
| 9 | Strives to promptly correct any errors discovered while producing the final report or after publication. As appropriate, disseminates the correction publicly or to others relying on the results. | 97.6 |
| 11 | Understands and conforms to confidentiality requirements of data collection, release, and dissemination and any restrictions on its use established by the data provider (to the extent legally required), protecting use and disclosure of data accordingly. | 97.6 |
| 14 | Is honest about their qualifications and about any limitations in their competence to complete a task. They accept full | 95.2 |
| | responsibility for their professional performance and practice. | |
| 16 | Should be forthright about any circumstances that might lead to either real or perceived conflicts of interest or otherwise tend to undermine the independence of their judgment. | 95.6 |
| 19 | Maintains high standards of professional competence, conduct, and ethical practices. | 96.5 (27) |
| 20 | Avoids plagiarism. The knowing presentation of another person's mathematical discovery as one's own constitutes plagiarism and is a serious violation of professional ethics. Plagiarism may occur for any type of work, whether written or oral and whether published or not. | 99.1 (27) |
| 22 | Avoids condoning or appearing to condone mathematical, scientific, or professional misconduct. | 92.6 (34) |
| 23 | Respects others; seeks and respects diverse opinions; promotes the equal dignity and fair treatment of all people; and neither engages in nor condones discrimination based on personal characteristics. Respects personal boundaries in interactions, and avoids harassment, including sexual harassment; bullying; and other abuses of power or authority. Takes appropriate action when aware of unethical practices by others. | 93.5 |
| 25 | Helps strengthen the work of others through appropriate peer review; in peer review, one assesses methods, not individuals. Strives to complete review assignments thoroughly, thoughtfully, and promptly. | 86.2 (33) |
| 26 | Avoids, and acts to discourage, retaliation against or damage to the employability of those who responsibly call attention to possible mathematical error or to scientific or other misconduct. | 97.2 (34) |
| 32 | Respects others and promotes justice, and inclusiveness in all work. Fosters fair participation of all people. Avoids and mitigates bias and prejudice. Does nothing to limit fair access. | 85.3 |
| 33 | Minimizes the possibility of indirectly or unintentionally harming others, mathematics practitioners should follow generally accepted best practices in academia, industry, and research, unless there is a compelling ethical reason to do otherwise. | 85.3 |
| 35 | Upholds, promotes, and respects the ethical responsibilities of the mathematics community. | 90.3 |
| 36 | Avoids and addresses exclusionary practices in hiring, teaching, and recruiting. | 92.1 (41) |
| | | |
| | An ethical mathematics practitioner who is a leader, employer, supervisor, mentor, or instructor... | |
| 41 | Maintains a working environment free from intimidation, including discrimination based on personal characteristics; bullying; coercion; unwelcome physical (including sexual) contact; and other forms of harassment. | 97.0 |
| 42 | Supports sound mathematical practice and exposes incompetent or corrupt mathematical practice. | 89.0 |
| 43 | Strives to protect the professional freedom and responsibility of mathematical practitioners who comply with these guidelines. | 90.1 (41) |
| 47 | Ensures that they enhance, not degrade, the quality of working life. Leaders should consider accessibility, physical safety, psychological well-being, and human dignity of all community members. | 93.9 (44) |
| 48 | Articulates, applies, and supports policies and processes that reflect the principles of mathematicians' responsibilities. Designing or implementing policies that deliberately or negligently violate, or tend to enable the violation of, mathematicians' responsibilities is ethically unacceptable. | 90.8 |
| 51 | Takes full responsibility for their contributions to the certification/granting of a degree in mathematics by ensuring the high level and originality of the Ph.D. dissertation work, and sufficient knowledge in the recipient of important branches of mathematics outside the scope of the thesis. | 86.7 (44) |

Table 5 presents endorsement results for the 17 items with levels of agreement between 70% and 84.9%.

Table 5. **Survey items with 70%-84.9% agreement**

| | The ethical mathematics practitioner... | % agree (# not responding) |
|---|---|---|
| 5 | Strives to follow, and encourages all collaborators to follow, an established protocol for authorship. | 83.7 |
| 8 | Assesses, and is transparent about, the origin and source of the tools and methods they use, including prior results and | 84.8 (17) |
| | data. Practitioners, when possible, acknowledge and disclose the origin of the problems they are solving and the interests that their work is intended to serve. | |
| 12 | Strives to ensure that data sources, choice of methods, and applications do not create or perpetuate social biases or discrimination. Seeks to avoid confirmation bias. | 84.8 |
| 13 | Recognizes any mathematical descriptions of groups may carry risks of stereotypes and stigmatization. Practitioners should contemplate, and be sensitive to, the manner in which information in their work across education, research, public policy, and in the public in general, is framed to avoid disproportionate harm to vulnerable groups | 80.0 |
| 18 | Strives to support and achieve quality work in both the process and products of professional work. | 82.5 |
| 21 | Understands the differences between questionable mathematical, scientific, or professional practices and practices that constitute misconduct. The ethical mathematics practitioner avoids all of the above and knows how each should be handled. | 81.6 (28) |
| 24 | Promotes sharing of data, methods, scholarship as much as possible and as appropriate without compromising propriety. | 79.8 (33) |
| 27 | Must know how to work ethically in collaborative environment by staying informed about and adhering to applicable rules, approvals, and guidelines to projects. Science and mathematical practice are often conducted in teams made up of practitioners with different professional and ethical standards. | 83.3 |
| 29 | When conducting their work in conjunction with other professions, must continue to abide by mathematicians' responsibilities, as well as any guidelines of the other professions. When there is a conflict or an absence in the partner profession's guidelines, the mathematical practitioners' responsibilities should be followed. | 78.6 (39) |
| 31 | Strives to resist institutional confirmation bias and systematic injustice. Opposes marginalization of people on the basis of human differences. When assessing or evaluating mathematics practitioners or their work, uses relevant subject matter-specific qualifications. Uses qualifications, performance, and contributions as the basis for decisions regarding mathematical practitioners of all levels. | 81.4 (40) |

| 39 | When involved in advising graduate students, should fully inform them about the employment prospects they may face upon completion of their degrees. | 80.4 (40) |
|---|---|---|
| | | |
| | An ethical mathematics practitioner who is a leader, employer, supervisor, mentor, or instructor... | |
| 40 | Recognizes that mathematicians' responsibilities exist and were articulated for the protection and support of the mathematics practitioner, the mathematics user, and the public alike. | 75.0 (42) |
| 44 | Recognizes the inclusion of mathematics practitioners as authors, or acknowledgement of their contributions to projects or publications, requires their explicit permission because it implies endorsement of the work. | 84.2 |
| 46 | Articulates and encourages acceptance and fulfillment of responsibilities by members of the organization or group. | 71.0 (42) |
| 49 | Ensures that opportunities are available to mathematics practitioners to help them improve their knowledge and skills in the practice and dissemination of mathematics, in the practice of ethics, and in their specific fields, and encourages people to take those opportunities. | 70.7 (43) |
| 50 | Demonstrates and educates students, employees, and peers on the ethical aspects of their teaching, ethical implications of their work, and the ethical challenges within the practice of mathematics. | 82.8 |
| 52 | Does not exploit the offer of a temporary position at an unreasonably low salary and/or an unreasonably heavy workload. | 78.6 |

Table 6 presents endorsement results for the 9 items with the lowest levels of agreement (38.8%-69.9%).

Table 6. **Survey items with 38.8%-69.9% agreement**

| | The ethical mathematics practitioner... | % agree (# not responding) |
|---|---|---|
| 10 | Strives to make new mathematical knowledge widely available to provide benefits to society at large and beyond their own scope of applications. | 57.6 |
| 15 | Recognizes that if they engage in mathematics practice, they do so in a social and cultural context. Practitioners | 56.1 (28) |
| | should contribute to society and to human well-being, acknowledging that all people are stakeholders in mathematics. | |
| 17 | Reviews submissions for peer review publication for potentially damaging/negative/unjust or inequitable implications. | 68.4 |
| 28 | Recognizes other professions have standards and obligations, research practices and standards can differ across disciplines, and statisticians (*sic*) do not have obligations to standards of other professions that conflict with mathematicians' responsibilities. | 51.9 |
| 30 | Instills in students and non-mathematicians an appreciation for the practical value of the concepts and methods they are learning or using. | 38.8 |
| 34 | Improves public awareness and understanding of mathematics and quantitative argument, related technologies, and their consequences. | 51.5 (39) |
| 37 | Values peoples' identity as part of their work. | 53.5 |
| 38 | Builds compassionate, sustaining community which is accountable to its members. Accepts their accountability to improve community. | 67.3 |
| | | |
| | An ethical mathematics practitioner who is a leader, employer, supervisor, mentor, or instructor... | |
| 45 | Encourages full participation of practicing mathematicians in recognizing their responsibilities, and encourages the recognition that one practices mathematics in a social context, not in value-free isolation. | 65.3 |

**Respondent suggestions for additional items**

In response to the open-ended item, “Please describe what you think is missing from the preceding list of items", we received 39 unique responses, several with multiple items listed. Of these 39 responses, nearly half (18) did not include specific items that respondents thought were missing (examples of these comments were “nothing is missing” and “You treat the word "Ethical" as if it (is) rigidly, naturally or easily defined. It is not”). Eleven of the respondents commented about clarity of items, or the perspective taken by the researchers (and did not suggest additional items). Eight respondents indicated additional elements were needed that were relevant to teaching/education specifically. Six responses were not interpretable (e.g, “connection to related professional societies’ ethical standards”), three were specific to the workplace (irrespective of what the work is), two (which included many suggestions) related to academic work (not teaching), including scholarship; two related to updating the ethical guidelines over time, and one referred to comments they made on other items earlier in the survey, but did not suggest anything was missing from the survey.

Informal thematic analysis of the responses to this open-ended item led to the following six general categories that represent domains/aspects beyond what might be considered to be the more abstract or objective topics of some mathematical work:

a. Workplace (not teaching, even if you do teach at work) – basic human respect/rights non-violations

b. Educating (if this is your primary job/role or if you teach only as part of mentoring/collaborating): teaching effectively, grading objectively; doing your best to promote learning.

c. Scholarship (writing, reviewing, and correcting errors, even if you made them) – respect for other’s work and others’ input to your work.

d. Respect for the profession/stewardship (apart from scholarship, ethical & objective reviewing)

e. Math in the world: effectively preparing learners & users of math; not gatekeeping.

f. Recognizing and effectively/respectfully treating stakeholders in work, teaching, scholarship, use of math, and the profession.

**Analysis of respondent comments on individual items**

Respondents were invited to comment on their answer to (or on the wording of) each item. The table below summarizes the themes derived from an informal analysis of the item-level comments. We treated comment content as a "theme" if it was observed in two or more comments on at least two items. The first column in Table 7 contains the comment theme we uncovered, the second gives the items in which the comment was made with the endorsement of that item in the first set of parentheses. In a second parenthetical next to each item in the second column, we give the number of comments coded with the theme out of the total number of comments for the given item.

Table 7. **Thematic analysis of comments on items, listed in order of frequency**

| **Comment Theme with Notes** | Item number (percent "YES")<br>(number of coded comments/total comments) |
| --- | --- |
| Certain terms are too vague/unclear/loaded/subjective<br><br>**Notes: Comment observed for items with a wide range of endorsement rates (and origins). Some readers want algorithmic guidelines, while ethical guidelines typically need room for interpretation in order to be comprehensive.** | 19 (96.5%) (3/7)<br>1 (94.3%) (6/17)<br>47 (93.9%) (3/9)<br>22 (92.6%) (4/10)<br>48 (90.8%) (3/7)<br>35 (90.3%) (3/11)<br>42 (89.0%) (2/15)<br>51 (86.7%) (4/15)<br>32 (85.3%) (4/18)<br>33 (85.3%) (4/13)<br>27 (83.3%) (2/14)<br>18 (82.5%) (5/10)<br>21 (81.6%) (5/27) |

| | 31 (81.4%) (4/20)<br>29 (78.6%) (2/19)<br>52 (78.6%) (11/27)<br>40 (75.0%) (3/14)<br>46 (71.0%) (10/16)<br>38 (67.3%) (2/14)<br>45 (65.0%) (4/18)<br>37 (53.5%) (8/29)<br>15 (56.1%) (2/34)<br>**22 total items** |
|---|---|
| Item is desirable but not an ethical obligation<br><br>**Notes: Comment observed for items with a wide range of endorsement rates (and origins). Further, some of these comments indicated that the participant thought the standard in the item was too high, and that the same item but with a lower standard might be an ethical obligation.** | 1 (94.3%) (2/17)<br>43 (91.0%) (2/9)<br>35 (90.3%) (3/11)<br>42 (89.0%) (3/15)<br>32 (85.3%) (7/18)<br>8 (84.8%) (4/20)<br>5 (83.7%) (3/27)<br>50 (82.8%) (4/12)<br>18 (82.5%) (4/10)<br>24 (79.8%) (7/16)<br>49 (70.7%) (3/15)<br>38 (67.3%) (3/14)<br>10 (57.6%) (8/31)<br>15 (56.1%) (2/34)<br>37 (53.5%) (2/29)<br>34 (51.5%) (21/28)<br>30 (38.8%) (16/35)<br>**17 total items** |
| Item does not apply to pure mathematics<br><br>**Notes: Comment observed for items with a wide range of endorsement rates (and origins).** | 7 (97.9%) (5/14)<br>11 (97.6%) (2/5)<br>16 (95.6%) (2/4)<br>3 (94.3%) (7/23)<br>8 (84.8%) (3/20)<br>12 (84.8%) (4/20)<br>5 (83.7%) (3/27)<br>17 (68.4%) (5/30)<br>10 (57.6%) (5/31)<br>37 (53.5%) (12/29)<br>**10 total items** |
| Standard is too difficult to meet (includes comments about external obstacles such as proprietary work as well as comments about standards that are intrinsically too difficult to meet such as identifying your own biases) | 20 (99.1%) (2/10)<br>7 (97.9%) (4/14)<br>9 (97.6%) (2/8)<br>3 (94.3%) (2/23) |

| **Notes: All but 1 of these pertains to items with 80% or higher support** | 8 (84.8%) (5/20)<br>12 (84.8%) (4/20)<br>21 (81.6%) (13/27)<br>10 (57.6%) (10/31)<br>**8 total items** |
|---|---|
| Item could harm those that should be protected (such as vulnerable or minoritized faculty)<br><br>**Notes: All of these pertain to items with 80% or higher support.** | 14 (95.2%) (3/8)<br>1 (94.3%) (2/17)<br>23 (93.5%) (3/16)<br>36 (92.1%) (4/11)<br>42 (89.0%) (5/15)<br>2 (87.9%) (2/28)<br>31 (81.4%) (3/20)<br>39 (80.4%) (2/24)<br>**8 total items** |
| Item contains 2 or more distinct guidelines (which may or may not contradict one another)<br><br>**Notes: Several items from the source codes have multiple statements in a single guideline. We chose not to separate them.** | 3 (94.3%) (2/23)<br>23 (93.5%) (2/16)<br>2 (87.9%) (10/28)<br>31 (81.4%) (9/20)<br>**4 total items** |
| Item is about managerial work, not mathematical work.<br><br>**Notes: Each of the four items where this theme emerged was about mathematicians in leadership roles.** | 47 (93.9%) (5/9)<br>36 (92.1%) (2/11)<br>52 (78.6%) (8/27)<br>49 (70.7%) (4/15)<br>**4 total items** |
| The truth is more important than the potential for harm.<br><br>**Notes: These comments may represent contemporary debates about speech and harm (the topics of the items). The items all have less than 85% endorsement support.** | 12 (84.8%) (5/20)<br>13 (80.0%) (5/20)<br>17 (68.4%) (5/30)<br>**3 total items** |
| Pure mathematics is not itself unethical/math is (can be) value-free<br><br>**Notes: All three of these items have endorsement rates of 65% or lower. There were a lot of comments(18-35) on each item.** | 45 (65.0%) (8/18)<br>15 (56.1%) (16/34)<br>30 (38.8%) (11/35)<br>**3 total items** |
| Item is tautological/leading/suggests its own answer<br><br>**Notes: Both items giving rise to this theme have endorsement rates above 90%.** | 20 (99.1%) (3/10)<br>35 (90.3%) (3/11)<br>**2 total items** |
| Wording of item is poor<br><br>**Notes: These comments refer to errors in the writing of the survey. In particular, for item 28, the authors mistakenly left the word "statistician" in the wording. Item 29 has the same essential meaning as item 28, but received 78.6% endorsement, whereas item 28 received 51.9% endorsement -and all 45 comments related to poor wording. These items may need to be rewritten or expanded for proto-guidelines.** | 26 (97.2%) (2/6)<br>28 (51.9%) (45/45)<br>**2 total items** |
| Item is not unique to mathematics<br><br>**Notes: There were several items in which this comment was made once.** | 41 (97.0%) (3/8)<br>23 (93.5%) (4/16)<br>**2 total items** |
| Item does not apply to all mathematicians<br><br>**Notes: One of these is along the lines of "managerial not mathematics" and the other is about public communication.** | 50 (82.8%) (2/12)<br>34 (51.5%) (7/28)<br>**2 total items** |

**Survey response drop-off**

Informal examination of missing data patterns suggested that if an item had a missing response, then all items following the first missing response were also missing. Quantifying the pattern showed that there was 0-1 response missing from the first few items, then 17 missing responses, with missingness jumping to 34, then 44. After considering the comments on individual items, and responses to the final open-ended item, "Please describe what you think is missing from the preceding list of items", we concluded that if respondents felt the first/first few items were unclear, or not relevant, or presented another challenge to responding, then that respondent would stop answering the survey items.

## Discussion

Our project was devised in order to answer three key questions to move the field forward. Our survey yielded the following results:

**1. Which elements of the existing ethical codes of the AMS and MAA are perceived to be relevant to ethical practice by the mathematics community?** Mathematics community members endorsed all of the 13 individual elements reflected on the AMS Code of Ethics that we included in our survey, 11 of which overlapped with items also included in the codes/guidelines of math-adjacent societies (ACM & ASA), yielding just three items unique to the AMS. One of the AMS items that is included by neither ACM nor ASA, "Does not exploit the offer of a temporary position at an unreasonably low salary and/or an unreasonably heavy workload", was endorsed by over 78% of respondents. As noted, our thematic analyses of the AMS and MAA guidelines led us to omit 100% of MAA and over 50% of AMS items from our survey. The rationale for these omissions was that the existing guideline elements were not relevant to "ethical mathematical practice" because they were too specific (e.g., to MAA employees or for AMS editorial roles). One reason for the specificity of the MAA and AMS guidelines - which led to these omissions - might be that they reflect a larger community practice and belief that mathematical practice is inherently neutral and value-free; thus any ethical guideline would have to do less with mathematical practice and more with specific roles (like MAA employment or AMS editorial duties). We saw evidence of these two perspectives in comments on the survey, echoing current scholarship (e.g., Buell & Piercey 2022, Pearson 2019, Shulman 2002).

**2. What ethical guidelines from math-adjacent professional societies (ACM, ASA) does the mathematics community believe are relevant?**

We adapted 42 of 52 items in this survey from the ACM and ASA ethical practice standards. Since we had to do a thematic analysis of MAA and AMS codes to create any items for the survey, we opted to use the more specific language of the ACM and ASA practice standards, except for the definition of plagiarism which was most concrete in the AMS code. The process by which we created the survey suggests that future ethical practice standards for mathematics should utilize a more elemental, possibly less narrative, approach to guidelines. Importantly, every item on the survey was recognizable to respondents - to the extent that they were able to either endorse it, reject it as relevant to ethical mathematical practice, comment on it, or some combination of these. We interpret this to mean that a typical practitioner would be able to find specific guidance, and possibly justify a course of action, with a more elemental representation of their ethical obligations in any given case or aspect of practice. So, both specific content and the organization of math-adjacent guidelines are relevant for the community.

Of the 52 items included, 51 of them were endorsed by more than 50% of the sample. The sole item that the majority did *not* endorse (only 38.8% endorsed it) was, “Instills in students and non-mathematicians an appreciation for the practical value of the concepts and methods they are learning or using.” This item was taken from the ASA Guidelines and modified for mathematics, but we neglected to consider the role of “practical value” in the way mathematical concepts are viewed. Eleven of the 35 comments in the responses to this item specifically expressed objections to this term.

One commenter even compared mathematics work to arts and humanities. While the idea of “practical value” has importance for statistics instruction, it has a different interpretation for mathematical practice and instruction. More generally, common reasons for the rejection of items by individual participants, based on our analysis of item-by-item comments, included the vagueness of terms, concern that an item was desirable but did not constitute “an ethical obligation”, and the perception that the item does not apply to “pure” mathematicians.

The ACM and ASA practice standards have been developed and refined by groups with the sole purpose of contemplating the relevance, wording, and applicability of ethical guidelines to both the practitioner and the practice itself. We asked respondents to consider both what they do, and the profession itself, in their consideration of whether an item is relevant to “ethical mathematical practice”. Going forward, these results suggest that useful input from these standards in terms of both content and the organization of elements (rather than narrative) can fruitfully be leveraged in the development of new guidelines for ethical mathematical practice.

In our discussions of the thematic analyses of comments, we determined that at least some of the reasons for those in the minority that voted “no” on an item reflected a need for balance in the drafting of guidelines: precision vs. flexibility in terms. Thus, finding the right minimum standard for what is “ethical” will require focused discussion. Other reasons, such as “item does not apply to pure mathematics” or “the item is managerial and not mathematical”, relate to considerations of professional identity and the role of ethical guidelines for a profession rather than an individual.

**3. What other guidelines are necessary that are unique to mathematics? What ethical guidance is lacking from the AMS, MAA, ACM, and ASA guidelines?**

Only one item had 100% agreement: "Discloses conflicts of interest, financial and otherwise, and manages or resolves them according to established (institutional/regional/local) rules and laws." This was included in some form on all four source documents. Our analysis of the endorsement rates and comments offered on each of the items suggested that there are important aspects of math-adjacent professions currently missing from existing guidance for ethical mathematical practice, but also highlighted important differences between disciplinary perspectives. For example, comments on 23 of the 52 items reflected a need for greater precision of language, and potentially less opportunity for subjectivity in the articulation of ethical obligations. Comments on 17 of the 52 items suggested that the perception of an "ethical obligation" may differ slightly for mathematics as compared to statistics or computing. Comments on 10 of 52 items suggests that, for at least some respondents, there is a distinction between the ethical obligations incurred in "pure" mathematics and those incurred in other types of mathematics. In terms of "what was missing", we noted six themes arising from the 39 suggestions for what was missing from our 52 items: workplace; teaching/grading/mentoring; scholarship; professional respect; effective preparation of users of mathematics (who are not mathematicians); and respect for stakeholders.

Important limitations to this study must be noted, chief of which is that we did not have a random sample of responses to our survey. However, the origins of the codes of ethics (MAA, AMS, ACM) and ethical guidelines (ASA) also tend to originate from a small cohort of individuals charged specifically with the task of creating or revising/maintaining the ethical practice standard or code. That is, while the ACM (2018) and ASA (2018/2022) explicitly sought input on the guidance documents from members of their respective organizations, none of the organizations (has) ever conducted a survey to gauge community endorsement of their ethical practice standards like the one we created and deployed. Instead, these organizations selected a small group of experienced practitioners and tasked them with generating, and/or revising, the guidance. Even with the sample of respondents we collected, we were able to answer our three research questions, as well as generating evidence that a mathematics-specific code of ethics can comprise more than the topics that have been included to date (e.g., AMS/MAA).

**Conclusions**

This study described the first national survey of mathematics practitioners to inquire about perceived relevance of elements of the AMS ethical code together with elements taken from math-adjacent professional societies for computing (ACM 2018) and statistics (ASA 2018). Our independent evaluation of the MAA and AMS codes yielded three unique-to-AMS items, and zero MAA items that were either not also included in other codes, or were deemed too specific to general mathematical practice to be included in the survey intended to outline ethical mathematical practice. Of the three items derived from the AMS, all were endorsed by more than half of our sample. Another 38 items, including four contributed from our JMM Town Hall attendees, and others adapted from the ASA and ACM plus one derived from both ASA & ACM, were also endorsed by more than half our sample with one exception: one item from the ASA

Ethical Guidelines was endorsed by less than 40% of the mathematicians surveyed. Although 51 (of 52) items were endorsed by the majority (over 50%) of respondents, 20/52 (38%) were endorsed by 90-100% of respondents, and 43/52 (82.7%) were endorsed by at least 70% of respondents. Nine of 52 (17.3%) items were endorsed by 38.8-69.9% of respondents, representing the clearest targets for further discussion.

When we examined suggestions for “what else” should be included in a comprehensive representation of “ethical mathematical practice”, the only specific suggestions that were different from what was already presented related to an ethical obligation to teach mathematics effectively. Although some comments suggested that abstract aspects of mathematics may complicate the relevance of ethical practice standards, nothing unique to mathematics or the practice/teaching of mathematics was suggested.

Our survey was intended to answer specific questions about perceptions of relevance for ethical guidance from AMS, MAA, and ethical practice standards from adjacent disciplines, and can only be viewed as a first step in the effort in creating ethical practice standards for mathematics. This is the topic of our ongoing work.

Acknowledgement:
This work was funded by a Collaborative Incubation Grant from NSF to the co-authors. We were humbled by the thoughtful, meaningful, and impactful comments from the community. They provided a foundation for a strong understanding of the practice and the possible next step regarding ethical guidelines for the practice.

**References**

Association for Computing Machinery (ACM). (2018) *Code of Ethics*, downloaded from https://www.acm.org/about-acm/code-of-ethics on 12 October 2018.

American Mathematical Society (AMS), (2005; updated 2019). *Ethical Guidelines of the American Mathematical Society*. Retrieved Feb. 8, 2019 from: https://www.ams.org/about-us/governance/policy-statements/sec-ethics

American Mathematical Society (AMS), (2007). *American Mathematical Society Policy on Supportive Practices and Ethics in the Employment of Young People in the Mathematical Sciences*. Retrieved Feb. 8, 2019 from: https://www.ams.org/about-us/governance/policy-statements/sec-supportivepractices

American Statistical Association (ASA), (2018; revised 2022). *ASA Ethical Guidelines for Statistical Practice-revised, Retrieved* 30 April 2018 from https://www.amstat.org/ASA/Your-Career/Ethical-Guidelines-for-Statistical-Practice.aspx

Antes AL, Wang X, Mumford MD, Brown RP, Connelly S, Devenport LD. (2010). Evaluating the effects that existing instruction on responsible conduct of research has on ethical decision making. *Academic Medicine* 85: 519-526.

Atweh, B., Vale, C., and Walshaw, M. (2012). Equity, diversity, social justice and ethics in mathematics education, Research in Mathematics Education in Australasia 2008-2011, Sense Publishers: Rotterdam, The Netherlands, pp.39-65.

Bass H. (2006). Developing scholars and professionals: the case of mathematics. In CM Golde, GE Walker (Eds). *Envisioning the Future of Doctoral Education: Preparing stewards of the discipline*. San Francisco: Jossey Bass.

Beauchamp, T.L. and N.E. Bowie, eds. (1979). *Ethical Theory and Business*, Prentice-Hall, Inc.: Englewood Cliffs, NJ.

Briggle A, & Mitcham C. (2012). *Ethics and science: An introduction*. Cambridge, UK: Cambridge University Press.

Buell, C. A. and Piercey, V. (2019) Ethics in Mathematics. *FOCUS*, Feb/March.

Buell, C. & Piercey, V. (2022). Special Issue -- Ethics in Mathematics, Eds. Buell and Piercey, *Journal of Humanistic Mathematics*, Volume 12 Issue 2. Available at: https://scholarship.claremont.edu/jhm/vol12/iss2/3

Chiodo, M. and Bursill-Hall P. (2018). Four Levels of Ethical Engagement. *EiM Discussion Papers, 1*. Cambridge University. Retrieved from https://ethics.maths.cam.ac.uk/assets/dp/18_1.pdf Feb. 8, 2019.

Committee on the Undergraduate Program in Mathematics (CUPM). (2015). 2015 CUPM Curriculum Guide to Majors in the Mathematical Sciences. Mathematical Association of America: Washington DC. Downloaded from https://www.maa.org/sites/default/files/CUPM%20Guide.pdf Feb. 8, 2019.

Ernest, P. (2018). The ethics of mathematics: Is mathematics harmful? In Ernest, P. ed., *The Philosophy of Mathematics Education Today*, ICME013 Monographs, Springer International Publishing AG: Basel, Switzerland, pp. 187-216.

Ferrini-Mundy J. (2008). What core knowledge do doctoral students in mathematics education need to know? In RE Reys, JA Dossey (Eds). *U.S. Doctorates in Mathematics Education: Developing Stewards of the Discipline*. Providence, RI: American Mathematical Society. pp. 63-74.

Franklin, J. (2005). A “Professional issues and ethics in mathematics” course. *Gazette of the Australian Mathematical Society,* 32 (2005), 98–100.

Hoffman, W.M., Moore, J.M., and Fredo, D.A. eds. (1984). *Corporate Governance and Institutionalizing Ethics: Proceedings of the Fifth National Conference on Business Ethics*, Lexington Books: Lexington, MA.

Kalichman MW. (2013). Why teach research ethics? In, National Academy of Engineering (Eds). *Practical Guidance on Science and Engineering Ethics Education for Instructors and Administrators*. Washington, DC: National Academies Press. Pp. 5-16.

Karaali, G. Doing Math in Jest: Reflections on Useless Math, the Unreasonable Effectiveness of Mathematics, and the Ethical Obligations of Mathematicians. *Mathematical Intelligence*r, Volume 41 Issue 3 (September 2019), pages 10-13.

Mathematical Association of America (MAA), (2021). MAA Code of Conduct In Support of a Welcoming and Inclusive Community. Retrieved from https://www.maa.org/about-maa/policies-and-procedures/maa-code-of-conduct on 5/19/2022.

Mathematical Association of America (MAA), (2017a). Welcoming Environment, Code of Ethics, and Whistleblower Policy. Retrieved from https://www.maa.org/about-maa/policies-and-procedures/welcoming-environment-code-of-ethics-and-whistleblower-policy on Feb. 8, 2019.

Mathematical Association of America (MAA), (2017b). *Policy on Conflict of Interest*. Retrieved from https://www.maa.org/about-maa/policies-and-procedures/policy-on-conflict-of-interest on Feb. 8, 2019.

May DR, Luth MT. (2013). The effectiveness of ethics education: a quasi-experimental field study. *Science and Engineering Ethics* 19(2): 545-568.

McNamara A, Smith J, & Murphy-Hill E. (2018). Does ACM's Code of Ethics Change Ethical Decision Making in Software Development? In *Proceedings of the 26th ACM Joint European Software Engineering Conference and Symposium on the Foundations of Software Engineering (ESEC/FSE '18), November 4–9, 2018, Lake Buena Vista, FL, USA.* ACM, New York, NY, USA, 5 pages. Downloaded from https://doi.org/10.1145/3236024.3264833.

Müller, D. (2018). Is there ethics in pure mathematics? *EiM Discussion Papers, 2*. Cambridge University. Accessed from https://ethics.maths.cam.ac.uk/assets/dp/18_2.pdf Feb. 8, 2019.

Neyland, J. (2004). Towards a postmodern ethics of mathematics education. In Walshaw, M. ed. (2004) *Mathematics Education Within the Postmodern.* Information Age Publishing: Greenwich, CT, pp. 55–73.

Neyland, J. (2008). Globalisation, ethics and mathematics education. In Atweh B., Calabrese Barton, A., Borba, M.C., Gough, N., Keitel-Kreidt, C., Vistro-Yu, C., Vithal, R., eds. *Internationalisation and Globalisation in Mathematics and Science Education*, Springer: Dordrecht, The Netherlands, pp. 113-128.

O'Neil, C. (2016). *Weapons of Math Destruction: How Big Data Increases Inequality and Threatens Democracy*. Crown Publishing Group: New York, NY.

Pearson, M. (2019, Oct.8, Nov. 26, Dec. 27) Parts I-III: The critical study of ethics in mathematics. Math Values Blog. Part I: https://www.mathvalues.org/masterblog/part-i-the-critical-study-of-ethics-in-mathematics Part II: https://www.mathvalues.org/masterblog/part-ii-the-critical-study-of-ethics-in-mathematics Part III: https://www.mathvalues.org/masterblog/part-iii-the-critical-study-of-ethics-in-mathematics

Rios CM, Golde C, & Tractenberg RE. (2019). The preparation of stewards with the Mastery Rubric for Stewardship: Re-envisioning the formation of scholars *and* practitioners. *Education Sciences 9(4), 292;* ***https://doi.org/10.3390/educsci9040292*** *Original* published on SocArXiv 7 jan 2019, DOI: 10.31235/osf.io/vw7j5.

Rittberg CJ, Tanswell FS, & Van Bendegem JP. (2020). Epistemic injustice in mathematics.*Synthese* 197, 3875–3904. https://doi.org/10.1007/s11229-018-01981-1

Shulman B. (2002) Is There Enough Poison Gas to Kill the City?: The Teaching of Ethics in Mathematics Classes, The College Mathematics Journal, 33:2, 118-125, DOI: 10.1080/07468342.2002.11921929

Sowder J.T. (1998). Ethics in mathematics education research, In Sierpinska A. and Kilpatrick J., eds. *Mathematics Education as a Research Domain: A Search for Identity*, New ICMI Studies Series, Vol 4. Springer: Dordrecht, The Netherlands, pp. 427-442.

Tractenberg RE. (2022-A). *Ethical Reasoning for a Data Centered World*. Cambridge, UK: Ethics International Press.

Tractenberg RE. (2022-B). *Ethical Practice of Statistics and Data Science*. Cambridge, UK: Ethics International Press.

Tractenberg RE. (2019, April 23). *Preprint.* Teaching and learning about ethical practice: The case analysis. Published in the *Open Archive of the Social Sciences* (SocArXiv), https://doi.org/10.31235/osf.io/58umw

Tractenberg RE, Russell A, Morgan G, FitzGerald KT, Collmann J, Vinsel L, Steinmann M, Dolling LM. (2015) Amplifying the reach and resonance of ethical codes of conduct through ethical reasoning: preparation of Big Data users for professional practice. *Science and Engineering Ethics*. 21(6):1485-1507. http://link.springer.com/article/10.1007%2Fs11948-014-9613-1 PMID: 25431219

Tractenberg RE, Cao J, Weisman M, Gillikin J, Rotelli M. (2021). Results of the first 5-yearly revision by the American Statistical Association's Committee on Professional Ethics Working Group on Revisions. In *Proceedings of the 2021 Joint Statistical Meetings*, held virtually. Alexandria, VA: American Statistical Association.

Weller, S. (1988). The effectiveness of corporate codes of ethics. *Journal of Business Ethics*. Vol 7, Issue 5, pp. 389 - 395.

## APPENDIX. Full survey results with item origins

Table A presents the results (percent of sample endorsing each item) for the 52 items on our survey, together with the origins of each item (and comments about changes, where relevant or noteworthy).

Table A

| | The ethical mathematics practitioner... | % agree (# not responding) |
|---|---|---|
| 1 | Works in a manner intended to produce valid, interpretable, and reproducible results. | 94.3 (1) |
| | **Source: ASA 2018** | |
| 2 | Does not knowingly accept work for which they are not sufficiently qualified, is honest about any limitation of expertise, and consults others when necessary or in doubt. It is essential that mathematics practitioners treat others with respect. | 87.9 |
| | **Source: ASA 2018** (ACM 2018 “Maintain high standards of professional competence, conduct, and ethical practice.” And AMS 2019 “Be knowledgeable in your field”) | |
| 3 | Identifies and mitigates any efforts to predetermine or influence the results or outcomes of mathematical practices; resists pressure to solve unethical problems/support predetermined outcomes. | 94.3 |
| | **Source: ASA 2018** | |
| 4 | Accepts full responsibility for their own work; does not take credit for the work of others; and gives credit to those who contribute. Respects and acknowledges the intellectual property of others. | 99.3 |
| | **Source: ASA 2018;** ACM 2018 (“respect the work required to produce new ideas…”); AMS 2019 (“Do not plagiarize…”) | |
| 5 | Strives to follow, and encourages all collaborators to follow, an established protocol for authorship. | 83.7 |
| | **Source: ASA 2018; AMS 2019** | |
| 6 | Discloses conflicts of interest, financial and otherwise, and manages or resolves them according to established (institutional/regional/local) rules and laws. | 100 |
| | **Source: ASA 2018; AMS 2019** | |
| 7 | Is candid about any known or suspected limitations, assumptions, or biases when working with data, methods, or models. Objective and valid interpretation of the results requires that the underlying analysis recognizes and acknowledges the degree of reliability and integrity of the data or the model. | 97.9 |
| | **Source: ASA 2018** | |
| 8 | Assesses, and is transparent about, the origin and source of the tools and methods they use, including prior results and data. Practitioners, when possible, acknowledge and disclose the origin of the problems they are solving and the interests that their work is intended to serve. | 84.8 (17) |
| | **Source: ASA 2018** | |
| 9 | Strives to promptly correct any errors discovered while producing the final report or after publication. As appropriate, disseminates the correction publicly or to others relying on the results. | 97.6 |
| | **Source: ASA 2018; AMS 2019** | |
| 10 | Strives to make new mathematical knowledge widely available to provide benefits to society at large and beyond their own scope of applications. | 57.6 |
| | **Source: ASA 2018** (AMS 2019) | |
| 11 | Understands and conforms to confidentiality requirements of data collection, release, and dissemination and any restrictions on its use established by the data provider (to the extent legally required), protecting use and disclosure of data accordingly. | 97.6 |
| | **Source: ASA 2018** (ACM 2018 “honor confidentiality”) | |
| 12 | Strives to ensure that data sources, choice of methods, and applications do not create or perpetuate social biases or discrimination. Seeks to avoid confirmation bias. | 84.8 |
| | **Source: ASA 2018** | |
| 13 | Recognizes any mathematical descriptions of groups may carry risks of stereotypes and stigmatization. Practitioners should contemplate, and be sensitive to, the manner in which information in their work across education, research, public policy, and in the public in general, is framed to avoid disproportionate harm to vulnerable groups | 80.0 |
| | **Source: ASA 2018** | |
| 14 | Is honest about their qualifications and about any limitations in their competence to complete a task. They accept full responsibility for their professional performance and practice. | 95.2 |
| | **Source: ACM 2018 (“is honest…”) and ASA 2018 (“they accept…”)** | |
| 15 | Recognizes that if they engage in mathematics practice, they do so in a social and cultural context. Practitioners should contribute to society and to human well-being, acknowledging that all people are stakeholders in mathematics. | 56.1 (28) |
| | **Source: ACM 2018** | |
| 16 | Should be forthright about any circumstances that might lead to either real or perceived conflicts of interest or otherwise tend to undermine the independence of their judgment. | 95.6 |
| | **Source: ACM 2018** (ASA 2018; AMS 2019) | |
| 17 | Reviews submissions for peer review publication for potentially damaging/negative/unjust or inequitable implications. | 68.4 |
| | **Source:** ACM (2018) (Thematic adaptation of 2.4 and 2.5) | |
| 18 | Strives to support and achieve quality work in both the process and products of professional work. | 82.5 |
| | **Source: ACM 2018** (ASA 2018) | |
| 19 | Maintains high standards of professional competence, conduct, and ethical practices. | 96.5 (27) |
| | **Source: ACM 2018** (ASA 2018; AMS 2019) | |
| 20 | Avoids plagiarism. The knowing presentation of another person's mathematical discovery as one's own constitutes plagiarism and is a serious violation of professional ethics. Plagiarism may occur for any type of work, whether written or oral and whether published or not. | 99.1 (27) |
| | **Source: AMS 2019** (ASA 2018; ACM 2018) | |
| | | |
| | The ethical mathematics practitioner... | |
| 21 | Understands the differences between questionable mathematical, scientific, or professional practices and practices that constitute misconduct. The ethical mathematics practitioner avoids all of the above and knows how each should be handled. | 81.6 (28) |
| | **Source: ASA 2018** | |
| 22 | Avoids condoning or appearing to condone mathematical, scientific, or professional misconduct. | 92.6 (34) |
| | **Source: ASA 2018** | |
| 23 | Respects others; seeks and respects diverse opinions; promotes the equal dignity and fair treatment of all people; and neither engages in nor condones discrimination based on personal characteristics. Respects personal boundaries in interactions, and avoids harassment, including sexual harassment; bullying; and other abuses of power or authority. Takes appropriate action when aware of unethical practices by others. | 93.5 |
| | **Source: ASA 2018** | |
| 24 | Promotes sharing of data, methods, scholarship as much as possible and as appropriate without compromising propriety. | 79.8 (33) |
| | **Source: ASA 2018** | |
| 25 | Helps strengthen the work of others through appropriate peer review; in peer review, one assesses methods, not individuals. Strives to complete review assignments thoroughly, thoughtfully, and promptly. | 86.2 (33) |
| | **Source: ASA 2018** (ACM 2018) | |
| 26 | Avoids, and acts to discourage, retaliation against or damage to the employability of those who responsibly call attention to possible mathematical error or to scientific or other misconduct. | 97.2 (34) |
| | **Source: ASA 2018** | |
| 27 | Must know how to work ethically in collaborative environment by staying informed about and adhering to applicable rules, approvals, and guidelines to projects. Science and mathematical practice are often conducted in teams made up of practitioners with different professional and ethical standards. | 83.3 |
| | **Source: ASA 2018** | |
| 28 | Recognizes other professions have standards and obligations, research practices and standards can differ across disciplines, and statisticians (sic) do not have obligations to standards of other professions that conflict with mathematicians' responsibilities. | 51.9 |
| | **Source: ASA 2018** | |
| 29 | When conducting their work in conjunction with other professions, must continue to abide by mathematicians' responsibilities, as well as any guidelines of the other professions. When there is a conflict or an absence in the partner profession's guidelines, the mathematical practitioners' responsibilities should be followed. | 78.6 (39) |
| | **Source: ASA 2018** | |
| 30 | Instills in students and non-mathematicians an appreciation for the practical value of the concepts and methods they are learning or using. | 38.8 |

| | **Source: ASA 2018** | |
|---|---|---|
| 31 | Strives to resist institutional confirmation bias and systematic injustice. Opposes marginalization of people on the basis of human differences. When assessing or evaluating mathematics practitioners or their work, uses relevant subject matter-specific qualifications. Uses qualifications, performance, and contributions as the basis for decisions regarding mathematical practitioners of all levels. | 81.4 (40) |
| | **Source: ASA 2018 "uses qualifications…" and "strives to resist confirmation bias") and JMM (community) 2020 "resist systematic injustice, opposes marginalizations…"** | |
| 32 | Respects others and promote justice, and inclusiveness in all work. Fosters fair participation of all people. Avoids and mitigates bias and prejudice. Does nothing to limit fair access. | 85.3 |
| | **Source: ACM 2018** (of the work, not related to education <where AMS 2019 has some similar language>) (ASA 2018, "respects others") | |
| 33 | Minimizes the possibility of indirectly or unintentionally harming others, mathematics practitioners should follow generally accepted best practices in academia, industry, and research, unless there is a compelling ethical reason to do otherwise. | 85.3 |
| | **Source: ACM 2018** | |
| 34 | Improves public awareness and understanding of mathematics and quantitative argument, related technologies, and their consequences. | 51.5 (39) |
| | **Source: ACM 2018** (ASA 2018) | |
| 35 | Upholds, promotes, and respects the ethical responsibilities of the mathematics community. | 90.3 |
| | **Source: ACM 2018** (ASA 2018) | |
| 36 | Avoids and addresses exclusionary practices in hiring, teaching, and recruiting. | 92.1 (41) |
| | **Source: JMM (community) 2020** | |
| 37 | Values peoples' identity as part of their work. | 53.5 |
| | **Source: JMM (community) 2020** | |
| 38 | Builds compassionate, sustaining community which is accountable to its members. Accepts their accountability to improve community. | 67.3 |
| | **Source: JMM (community) 2020** | |
| 39 | When involved in advising graduate students, should fully inform them about the employment prospects they may face upon completion of their degrees. | 80.4 (40) |
| | **Source: AMS 2019** | |
| | | |
| | An ethical mathematics practitioner who is a leader, employer, supervisor, mentor, or instructor... | |
| 40 | Recognizes that mathematicians' responsibilities exist and were articulated for the protection and support of the mathematics practitioner, the mathematics user, and the public alike. | 75.0 (42) |
| | **Source: ASA 2018** | |
| 41 | Maintains a working environment free from intimidation, including discrimination based on personal characteristics; bullying; coercion; unwelcome physical (including sexual) contact; and other forms of harassment. | 97.0 |
| | **Source: ASA 2018** (ACM 2018) | |
| 42 | Supports sound mathematical practice and exposes incompetent or corrupt mathematical practice. | 89.0 |
| | **Source: ASA 2018** | |
| 43 | Strives to protect the professional freedom and responsibility of mathematical practitioners who comply with these guidelines. | 90.1 (41) |
| | **Source: ASA 2018** | |
| 44 | Recognizes the inclusion of mathematics practitioners as authors, or acknowledgement of their contributions to projects or publications, requires their explicit permission because it implies endorsement of the work. | 84.2 |
| | **Source: ASA 2018** | |
| 45 | Encourages full participation of practicing mathematicians in recognizing their responsibilities, and encourages the recognition that one practices mathematics in a social context, not in value-free isolation. | 65.3 |
| | **Source: ACM 2018** | |
| 46 | Articulates and encourages acceptance and fulfillment of responsibilities by members of the organization or group. | 71.0 (42) |
| | **Source: ACM 2018** | |
| 47 | Ensures that they enhance, not degrade, the quality of working life. Leaders should consider accessibility, physical safety, psychological well-being, and human dignity of all community members. | 93.9 (44) |
| | **Source: ACM 2018** | |
| 48 | Articulates, applies, and supports policies and processes that reflect the principles of mathematicians' responsibilities. Designing or implementing policies that deliberately or negligently violate, or tend to enable the violation of, mathematicians' responsibilities is ethically unacceptable. | 90.8 |
| | **Source: ACM 2018** | |
| 49 | Ensures that opportunities are available to mathematics practitioners to help them improve their knowledge and skills in the practice and dissemination of mathematics, in the practice of ethics, and in their specific fields, and encourages people to take those opportunities. | 70.7 (43) |
| | **Source: ACM 2018** | |
| 50 | Demonstrates and educates students, employees, and peers on the ethical aspects of their teaching, ethical implications of their work, and the ethical challenges within the practice of mathematics. | 82.8 |
| | **Source: JMM (community) 2020** | |
| 51 | Takes full responsibility for their contributions to the certification/granting of a degree in mathematics by ensuring the high level and originality of the Ph.D. dissertation work, and sufficient knowledge in the recipient of important branches of mathematics outside the scope of the thesis. | 86.7 (44) |
| | **Source: AMS 2019** | |
| 52 | Does not exploit the offer of a temporary position at an unreasonably low salary and/or an unreasonably heavy workload. | 78.6 |
| | **Source: AMS 2019** | |